\documentclass[11pt]{article}  % Basic style 
\usepackage{amsthm}            % amsthm is used for theorem style (see below)
\usepackage{amssymb}           % AMS symbols are used
\usepackage{epsf}          % ESSENTIAL package for .eps (graphic) files  
\newcommand{\stdspace}{\hskip 0.75em plus 0.15em \ignorespaces} 
\newcommand{\ppar}{\par\goodbreak\vskip 8pt plus 3pt minus 3pt} 
\newcommand{\sh}[1]{\penalty-800\ppar{\bf #1}\par\medskip\nobreak}
\newcommand{\rk}[1]{\ppar{\bf #1}\stdspace}    
\makeatletter
\newtheoremstyle{plain}{14pt plus6.3pt minus6.3pt}{7.4pt plus3pt minus3pt}%
{\sl}{}{\bf}{}{0.75em}{\thmname{#1}\thmnumber{ #2}\thmnote{\rm\stdspace(#3)}}
\newtheoremstyle{definition}{14pt plus6.3pt minus6.3pt}{7.4pt plus3pt minus3pt}%
{\rm}{}{\bf}{}{0.75em}{\thmname{#1}\thmnumber{ #2}\thmnote{\sl\stdspace#3}}
\theoremstyle{plain}               
\renewenvironment{proof}[1][\proofname]{\par
  \normalfont
  \topsep\medskipamount \trivlist
  \item[\hskip\labelsep\bf
    #1]\ignorespaces
}{%
  \qed\endtrivlist\par
}
\def\sqr#1#2{{\vcenter{\vbox{\hrule  height.#2pt
	\hbox{\vrule width.#2pt height#1pt \kern#1pt \vrule width.#2pt}
	\hrule height.#2pt}}}}
\let\@footnote@\footnote
\def\footnote#1{\@footnote@{\small #1}}
\long\def\@caption#1[#2]#3{\par\addcontentsline{\csname
  ext@#1\endcsname}{#1}{\protect\numberline{\csname
  the#1\endcsname}{\ignorespaces #2}}\begingroup
    \@parboxrestore
    \small
    \@makecaption{\csname fnum@#1\endcsname}{\ignorespaces #3}\par
  \endgroup}
\let\@document@\document
\def\document{\@document@%
\setlength{\abovedisplayskip}{\medskipamount}
\setlength{\belowdisplayskip}{\medskipamount}}  
\let\@thebibliography@\thebibliography
\def\thebibliography#1 {\@thebibliography@{999}\small\parskip0pt % 
plus2pt\relax}
\let\@itemize@\itemize
\def\itemize{\@itemize@\parskip 0pt\relax}
\def\@listi{\leftmargin28.5pt\parsep 0pt\topsep 0pt 
 \itemsep4pt plus3pt minus2pt}
\let\@listI\@listi
\@listi
\makeatother
\def\relabelbox{%
  \hbox\bgroup%
}%
\def\endrelabelbox{%
\egroup%
}%
\def\relabel #1#2 {%
  \special{ps:/a {} def}%
  \smash{\rlap{#2}}%
}%
\def\adjustrelabel <#1,#2> #3#4 {%
  \special{ps:/a {} def}%
  \smash{\rlap{\kern #1 \raise #2\hbox{#4}}}%
}%
\def\extralabel <#1,#2> #3 {\smash{\rlap{\kern #1 \raise #2\hbox{#3}}}}%
\newcommand{\co}{\colon\thinspace}    %  Colon with good spacing for maps.
\newcommand{\np}{\newpage}            
\newcommand{\nl}{\hfil\break}         
\newcommand{\cl}{\centerline}         
\newcommand{\fnote}{\footnote}        
\newcommand{\inv}{^{-1}}              
\newcommand{\half}{{1 \over 2}}
\newcommand{\wrt}{with respect to }
\newcommand{\basl}{\,\backslash\,}
\newcommand{\up}{\uparrow} 
\newcommand{\down}{\downarrow}
\newcommand{\N}{{\mathbb N}} 
\newcommand{\Z}{{\mathbb Z}}
\newcommand{\R}{{\mathbb R}}

\newcommand{\C}{{\mathbb C}}
\newcommand{\G}{{\cal G}}
\newcommand{\pa}{\partial}

\newtheorem{lemma}{Lemma}[section]
\newtheorem{theorem}[lemma]{Theorem}
\newtheorem{prop}[lemma]{Proposition}
\newtheorem{cor}[lemma]{Corollary}
\newtheorem{lem}[lemma]{Lemma}

\newtheorem{algo}[lemma]{Algorithm}
\newtheorem*{Dthm}{Dehornoy's theorem}

\theoremstyle{remark}
\newtheorem{remark}[lemma]{Remark}

\def\title#1{\vglue 0.2truein\cl{\Large\bf #1}\vglue 0.15truein}
\def\author#1{{\parskip0pt\leftskip0pt plus 1fil\def\\{\par}{\sc#1}
\par\vglue0.1truein}}
\def\address#1{{\parskip0pt\leftskip0pt plus 1fil\def\\{\par}{\small#1}
\par\vglue7pt}}
\def\abstract{{\bf Abstract}\vglue 5pt}

\begin{document}

%   Title page 
%   (title, emails, abstract, class numbers and keywords
%   are all here, addresses are after the bibliography)
%
\title{Ordering the braid groups}

\author{Roger Fenn, Michael T~Greene, Dale Rolfsen\\
Colin Rourke, Bert Wiest}

\address{\rm Email:\stdspace \tt R.A.Fenn@sussex.ac.uk\ \ 
Michael.Greene@uk.radan.com\\
rolfsen@math.ubc.ca\ \ cpr@maths.warwick.ac.uk\\
bertw@gyptis.univ-mrs.fr}

\begin{abstract}

We give an explicit geometric argument that Artin's braid group $B_n$
is right-orderable. The construction is elementary, natural, and leads
to a new, effectively computable, canonical form for braids which we
call {\sl left-consistent canonical form}.  The left-consistent form
of a braid which is positive (respectively negative) in our order has
consistently positive (respectively negative) exponent in the smallest
braid generator which occurs.  It follows that our ordering is
identical to that of Dehornoy \cite{D.axiom}, constructed by very
different means, and we recover Dehornoy's main theorem that any braid
can be put into such a form using either positive or negative exponent
in the smallest generator but not both.

Our definition of order is strongly connected with Mosher's
normal form \cite{Mosher2} and this leads to an algorithm to decide 
whether a given braid is positive, trivial, or negative which is 
quadratic in the length of the braid word.

\end{abstract}
%
%  AMS classification numbers and keywords (hand-formatted):

{\bf AMS Classification numbers}\stdspace 
Primary:\stdspace 20F60, 06F15, 20F36\par
Secondary:\stdspace 57M07, 57M25
\vglue 5pt


\np  % page break at the end of the title page

\setcounter{section}{-1}
\section{Introduction}

Dehornoy \cite{D.proofnoaxiom,D.axiom,D.algorithm}
has proved that the braid group is right-orderable.  More precisely, 
there is a total order on the elements of
the braid group $B_n$ which is right invariant in the following sense.
Suppose that $\alpha$, $\beta$, $\gamma \in B_n$ and $\alpha < \beta$, then
$\alpha\gamma < \beta\gamma$.  This ordering is uniquely defined by
the condition that 
% $\sigma_i$ is infinitely greater than $\sigma_j$ if $j > i$.  
a braid $\beta_0\sigma_i \beta_1$ is positive (ie greater than the identity 
braid), where $\beta_0$, $\beta_1$ are words in $\sigma_{i+1}^{\pm 1},\ldots,
\sigma_{n-1}^{\pm 1}$. 
Dehornoy's proof is based on some highly complicated algebra
connected with left-distributive systems.  In this paper we construct
this order geometrically using elementary arguments.  

Our construction leads to a new, effectively 
computable, canonical form for braids which we call {\sl left-consistent
canonical form}.  The left-consistent form of a positive braid 
%(ie a braid greater than the identity braid) 
has the general shape  
$$
          \beta_0\sigma_i^e \beta_1\ldots \beta_{l-1} \sigma_i^e \beta_l
$$
where the $\beta_i$ are words in $\sigma_{i+1},\ldots,\sigma_{n-1}$ and their
inverses, and $e=+1$.  For a negative braid the form is similar but
with $e=-1$.  It follows at once that our ordering is identical to
Dehornoy's  and we recover Dehornoy's main theorem that any braid 
can be put into such a shape for $e=1$ or $e=-1$ but not both.

Our definition of order is strongly connected with Mosher's
automatic structure \cite{Mosher2} and this implies that the braid group
is {\sl order automatic}, ie the order can be detected from the
automatic normal form by a finite state automaton.  Furthermore the
resulting algorithm to decide whether a given braid is positive, 
trivial, or negative is linear in the length of the Mosher normal form 
and hence quadratic in the length of the braid word (in contrast, 
Dehornoy's algorithm \cite{D.algorithm}, although apparently fast 
in practice is only known to be exponential).

The paper is organised as follows.  Section 1 contains basic
definitions and introduces the curve diagram associated to a braid.
In section 2 we prove that a curve diagram can be placed in a unique
reduced form with respect to another and in section 3 we define the
order by comparing the two curve diagrams in reduced form, and prove
that it is right-invariant.  In section 4 we construct the
left-consistent canonical form of a braid, deduce that our order
coincides with Dehornoy's and recover Dehornoy's results. In section 5
we give some counterexamples connected with the order, and in section
6 we make the connection with Mosher's normal form and deduce the
existence of the quadratic time algorithm to detect order.  Finally,
in an appendix, we use {\sl cutting sequences} to give a formal
algorithm to turn a braid into the new left-consistent canonical form;
note that this algorithm is not quadratic time. 

\rk{Acknowledgements} We are grateful to the organisers of the
low-dimensional topology conference held at the Isle of Thorns in
Spring 1997, which was supported by the LMS, for providing a congenial
atmosphere for the initial work on this paper.  We are also grateful
to Caroline Series for suggesting that our curve diagrams might be
related to Mosher's normal form for mapping class groups.  We would
also like to thank the referee for helpful comments and a speedy
report.  Bert Wiest is supported by a TMR (Marie Curie) research
training grant.

%-------------------------------------------------------------------

\section{Braids and curve diagrams}

Let $D^2$ be the closed unit disk in $\C$, and let $D_n$ be the disk $D^2$
with $n$ distinct points in the real interval $(-1, 1)$ removed. We consider 
the group $B_n$ of self-homeomorphisms $\gamma\co D_n \to D_n$ with 
$\gamma|_{\partial D_n}=id$, up to isotopy of $D_n$ fixed on $\partial D_n$.
Multiplication in $B_n$ is defined by composition. 
The group $B_n$ is well-defined independently of the points removed;  indeed
if $D'_n$ is a disk with any $n$--tuple of points removed and $B'_n$ the
corresponding group then there is an isomorphism $B'_n\cong B_n$;
if these points also lie on $(-1,1)$ then this isomorphism is natural.

\begin{figure}[htb] %
\cl{
\relabelbox
\epsfbox{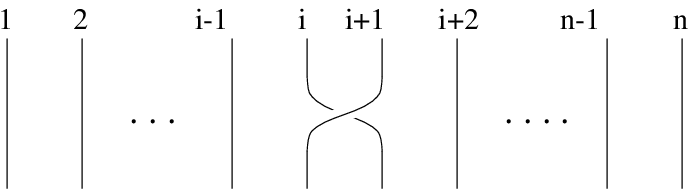} % exported at 60pc as portrait
\relabel{1}{\small $1$}
\relabel{2}{\small $2$}
\relabel{i-1}{\small $i-1$}
\relabel{i}{\small $i$}
\relabel{i+1}{\small $i+1$}
\relabel{i+2}{\small $i+2$}
\relabel{n-1}{\small $n-1$}
\relabel{n}{\small $n$}
\relabel{. . .}{$\ldots$}
\relabel{. . . .}{$\ldots$}
\endrelabelbox
}
\caption{The standard generator $\sigma_i$ of the braid group on $n$ strings}
\label{standgen}
\end{figure}

The group $B_n$ is isomorphic to the group $\widetilde{B}_n$ of braids on $n$
strings, with multiplication given by concatenation: if $\alpha,\beta$
are braids (pictured vertically) then $\alpha\cdot\beta$ is $\alpha$
above $\beta$. It is well known that this group has presentation
$$ 
\widetilde{B}_n\cong \langle \sigma_1,\ldots, \sigma_{n-1}\ | \ 
\sigma_i\sigma_j=\sigma_j\sigma_i \hbox{ \ if \ }|i-j| \geqslant 2,\ 
\sigma_{i+1}\sigma_i\sigma_{i+1}=\sigma_i\sigma_{i+1}\sigma_i \rangle,
$$
where the generator $\sigma_i$ ($i\in \{1,\ldots n-1\}$) is indicated in 
figure \ref{standgen}.

The isomorphism $\widetilde{B}_n \to B_n$ is given by `putting the braid 
in a solid cylinder and sliding $D_n$ once along it'. 
The inverse map is defined as follows: extend a given homeomorphism 
$\gamma\co D_n \to D_n$ to a homeomorphism $\gamma'\co D^2 \to D^2$,
then find a boundary-fixing isotopy $\gamma_t\co D^2 \to D^2$ with
$\gamma_0=id$ and $\gamma_1=\gamma$. Then the flow of the $n$ holes
of $D_n$ under $\gamma_t$ defines a braid on $n$ strings. 
For details see \cite{Birman}.

\begin{figure}[htb] 
\cl{
\relabelbox
\epsfbox{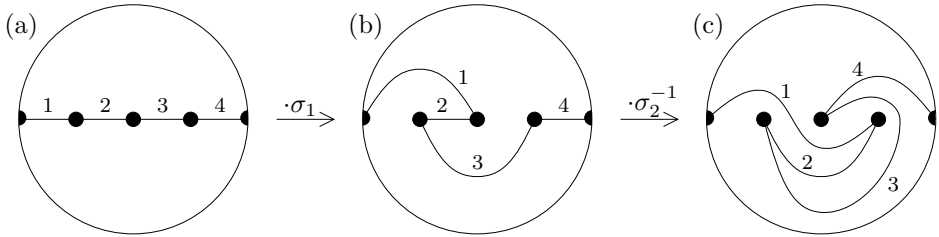} % exported at 60pc as portrait
\relabel{(a)}{\small (a)}
\relabel{(b)}{\small (b)}
\relabel{(c)}{\small (c)}
\relabel{0}{$\scriptstyle1$}
\relabel{1}{$\scriptstyle2$}
\relabel{2}{$\scriptstyle3$}
\relabel{3}{$\scriptstyle4$}
\relabel{0'}{$\scriptstyle1$}
\relabel{1'}{$\scriptstyle2$}
\relabel{2'}{$\scriptstyle3$}
\relabel{3'}{$\scriptstyle4$}
\relabel{0''}{$\scriptstyle1$}
\relabel{1''}{$\scriptstyle2$}
\relabel{2''}{$\scriptstyle3$}
\relabel{3''}{$\scriptstyle4$}
\adjustrelabel <-1pt,0pt> {xs_1}{\small $\cdot\sigma_1$}
\adjustrelabel <-1pt,1pt> {xs_2}{\small $\cdot\sigma_2\inv$}
\endrelabelbox
}
\caption{Curve diagrams of the braids $1$, $\sigma_1$, and 
$\sigma_1 \sigma_2\inv$}  \label{dicurv}
\end{figure}

On $D_n$ we draw $n+1$ line segments as in figure \ref{dicurv}(a). 
If $\gamma$ is a homeomorphism of $D_n$ representing an element $[\gamma]$
of $B_n$, then $\gamma$ sends these line segments to $n+1$ disjoint embedded
curves, and if $[\gamma_1]=[\gamma_2]\in B_n$ then $\gamma_1$ and $\gamma_2$ 
give rise to isotopic collections of curves. For instance, figure \ref{dicurv}
shows the effect of the braid $\sigma_1 \sigma_2\inv\in B_3$.
Here the holes of $D_n$, as well as $\pm 1$ are indicated by black dots.
We call such a diagram of $n+1$ disjoint simple curves in an $n$--punctured 
disk a {\sl curve diagram}, and we number the curves in the diagram $1$ to
$n+1$, as in figure \ref{dicurv}.

Conversely, from the curve diagram we can reconstruct the homeomorphism
$\gamma$ up to boundary-fixing isotopy; that is, we can reconstruct the
element of the braid group.

% ------------------------------------------------------------------------

\section{Reduced form} \label{setup}
Let $\Gamma$ and $\Delta$ be curve diagrams of two braids $\gamma$ and
$\delta$, say.  In order to compare $\Gamma$ and $\Delta$, we
superimpose the two diagrams and {\sl reduce} the situation by
removing unnecessary intersections.  This process is well known and
often called ``pulling tight'' (see eg \cite{Mosher2}).

We will denote the $i$th curve of a curve diagram such as $\Gamma$ by
$\Gamma_i$. The curves $\Gamma_i$ and the $\Delta_j$ are called
{\sl parallel} if they connect the same pairs of points and are isotopic 
in $D_n$. For instance, curve 3 of figure \ref{dicurv}(b) and curve 2 of
figure \ref{dicurv}(c) are parallel. We define $\Delta$ to be 
{\sl transverse} to $\Gamma$ if every curve of $\Delta$ 
either coincides precisely with some (parallel) curve of $\Gamma$,
or intersects the curves of $\Gamma$ transversely.
%(Note that the curve may or may not be parallel to a curve of $\Gamma$.)

We define the {\sl intersection index} of two transverse curve diagrams 
to be 
$$
n+1 + \#(\hbox{transverse intersections}) - \#(\hbox{coincident curves}).
$$
(The geometric significance is that $D_n$ cut along $\Gamma$ has two
components, and cutting {\it in addition} along $\Delta$ increases the
number of components by the intersection index.)  For example, the
diagrams in figure \ref{dicurv}(a) and \ref{dicurv}(c) have
intersection index 6, the diagrams in figure \ref{dicurv}(a) and
\ref{dicurv}(b) have intersection index $2$ and the diagrams in figure
\ref{dicurv}(b) and \ref{dicurv}(c) have intersection index $5$.  The
intersection index of two curve diagrams is $0$ if and only if the
diagrams are identical.

We now fix a curve-diagram $\Gamma$ for $\gamma$, and look at all possible
curve-diagrams for $\delta$. They are all isotopic in $D_n$, but they may 
have very different intersection-indices with $\Gamma$. We say $\Delta$ and 
$\Delta'$ are {\sl equivalent} (\wrt 
$\Gamma$) if they are related by an isotopy of $D_n$, which is fixed on 
curves of $\Delta$ which coincide with curves of $\Gamma$, and which leaves 
the diagrams transverse all the time. So coincident curves remain coincident, 
and the intersection index remains unchanged.

\begin{figure}[htb] 
\cl{
\relabelbox
\epsfbox{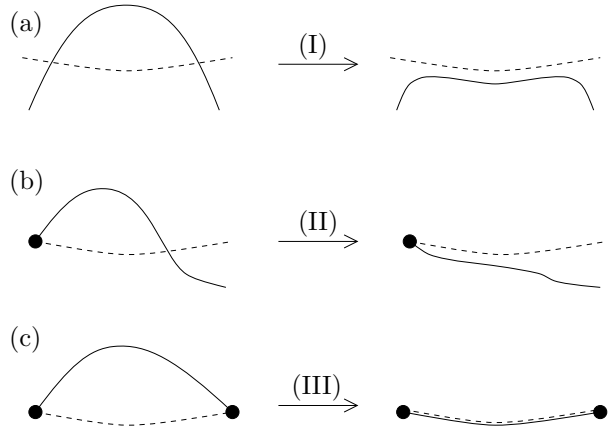} % exported at 55pc portrait
\relabel{(a)}{\small (a)}
\relabel{(b)}{\small (b)}
\relabel{(c)}{\small (c)}
\adjustrelabel <0pt,2pt> {(I)}{\small (I)}
\adjustrelabel <0pt,2pt> {(II)}{\small (II)}
\adjustrelabel <0pt,2pt> {(III)}{\small (III)}
\endrelabelbox
}
\caption{Three types of $D$--disks, and how to use them to reduce intersection
indices}
\label{Ddisk}
\end{figure}

We define a {\sl $D$--disk\fnote{$D$--disks are often called ``bigons''
in the literature.} between $\Delta$ and $\Gamma$} to be a subset of
$D_n$ homeomorphic to an open disk, which is bounded by one open
segment of some curve of $\Delta$, one open segment of some curve of
$\Gamma$, and two points, each of which may be an intersection-point
of the two curves or one of the `holes' of $D_n$, or $\pm 1\in D_n$.
%(We shall often leave out the words `\wrt $\Gamma$'.)  There are
three types of $D$--disks (types (a), (b), and (c)), indicated in
figure \ref{Ddisk}, where the curve-diagram $\Gamma$ is drawn with
dashed, and $\Delta$ with solid lines, and the dots denote holes or
$\pm 1$.

If there are no $D$--disks between $\Delta$ and $\Gamma$ then we say 
$\Delta$ and $\Gamma$ are {\sl reduced}. If the curve-diagrams $\Delta$ 
and $\Gamma$ are not reduced, ie if they have
a $D$--disk, then we can isotope $\Delta$ so as to reduce the 
intersection index with $\Gamma$ (figure \ref{Ddisk}). This isotopy consists 
of `sliding a segment of a curve of $\Delta$ across the $D$--disk' (for
reduction moves (I) and (II)), and of `squashing a $D$--disk to a line' (for 
reduction move (III)). The three moves reduce the intersection index by
2, 1, 1, respectively. We observe that any curve-diagram $\Delta$ with 
intersection index $0$ with $\Gamma$ is reduced. Thus we can reduce curve 
diagrams by a finite sequence of `isotopies across $D$--disks' as in 
figure \ref{Ddisk}.

\begin{lem}[Triple reduction lemma]\label{triplem}
Suppose $\Sigma$, $\Gamma$ and $\Delta$ are three curve diagrams such that
$\Gamma$ and $\Delta$ are both reduced \wrt $\Sigma$. Then there exists an
isotopy between $\Delta$ and a curve diagram $\Delta'$, which is an 
equivalence \wrt $\Sigma$, such that $\Sigma$, $\Gamma$ and $\Delta'$
are pairwise reduced.
\end{lem}

\begin{figure}[htb] 
\cl{
\relabelbox
\epsfbox{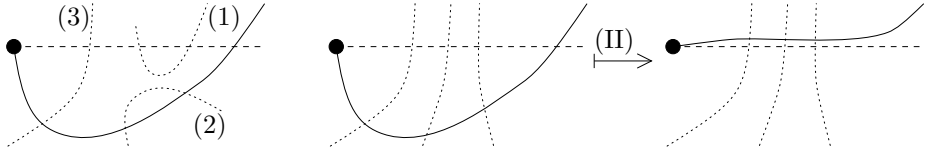}$$ 
\relabel{(1)}{\small (1)}
\relabel{(2)}{\small (2)}
\relabel{(3)}{\small (3)}
\adjustrelabel <0pt,2pt> {(II)}{\small (II)}
\endrelabelbox
}
% exported at 60pc portrait
\caption{The solid line is $\Delta$, the dashed $\Gamma$, and the dotted 
$\Sigma$} \label{tripred}
\end{figure}

\begin{proof} We consider a $D$--disk bounded by one segment of curve of
$\Gamma$ and one of $\Delta$. This $D$--disk may have several
intersections with $\Sigma$. There are, a priori, three possibilities
for the type of such an intersection --- they are indicated in figure
\ref{tripred}, labelled (1), (2), and (3). (In this figure, the
$D$--disk is of type (b), the cases of types (a) and (c) are similar.)

However, (1) and (2) are impossible, because $\Gamma$ and 
$\Delta$ are reduced \wrt $\Sigma$. So all intersections are of type
(3), and the $D$--disk can be removed without disturbing
the reduction of $\Sigma$ with respect to $\Gamma$ or $\Delta$, as indicated
in figure \ref{tripred}. The statement follows inductively. \end{proof}

\begin{lem}\label{reducedcoincide}
If two isotopic curve diagrams $\Gamma$ and $\Delta$ are reduced
\wrt each other, then they coincide.\end{lem}

\begin{proof} Suppose that the first curve $\Gamma_1$ of $\Gamma$ does not
coincide with the first curve $\Delta_1$ of $\Delta$.  Consider the
word obtained by reading the intersections of $\Gamma_1$ with the curves
of $\Delta$ in order.  Since $\Gamma_1$ is isotopic to $\Delta_1$, this word 
must cancel to the trivial word.  It follows by a simple innermost disk 
argument that there must be a $D$--disk. Hence $\Gamma_1$ must coincide with 
$\Delta_1$.  Similarly all curves of $\Gamma$ and $\Delta$ must coincide. \end{proof}

\begin{prop} \label{welldef} \sl 
If two curve diagrams $\Delta$ and $\Delta'$ of a braid $\delta$
are reduced \wrt $\Gamma$ then they are equivalent \wrt $\Gamma$.
\end{prop}

\begin{proof} By the triple reduction lemma we may reduce $\Delta$ \wrt
$\Delta'$ by an isotopy of $\Delta$ which is an equivalence \wrt
$\Gamma$.  After this reduction $\Delta$ and $\Delta'$ coincide by
lemma \ref{reducedcoincide}.  \end{proof}

We have proved that by reducing a curve diagram $\Delta$ \wrt a curve
diagram $\Gamma$ we can bring $\Delta$ into a uniquely defined
standard form \wrt $\Gamma$.  In particular reduction of $\Delta$ \wrt
the trivial curve diagram representing $1\in B_n$ (figure
\ref{dicurv}(a)) leads to a canonical representation of braids in
terms of {\sl cutting sequences}, which will be discussed in detail in
the appendix.

\begin{remark} \label{onecurv} \rm
The following observation will be crucial at a later point.
Let $\Gamma$ and $\Delta$ be transverse curve diagrams. Suppose the curve 
$\Delta_i$ on its own is reduced \wrt $\Gamma$.
%be the diagram consisting only of the $i$th curve of $\Delta$. 
%Suppose $\Delta_i$ is reduced \wrt $\Gamma$. 
Then we can reduce $\Delta$ \wrt $\Gamma$ by an isotopy of $\Delta$ which 
is fixed on $\Delta_i$. 
\end{remark}

% ------------------------------------------------------------------------

\section{The right-invariant order on $B_n$}
We define a total ordering on the braid group $B_n$ as follows.
Suppose $\gamma$ and $\delta$ are two braids on $n$ strings. We let $\Gamma$ 
be a curve diagram for $\gamma$, and $\Delta$ be a curve diagram for $\delta$
which is reduced \wrt $\Gamma$. The collection of curves of $\Gamma$ cuts
$D_n$ into two components, which we call the {\sl upper} and the {\sl lower}
component, containing the points $\sqrt{-1}$ respectively $-\sqrt{-1}$ in 
$D_n \subseteq \C$.
We orient the curves of $\Delta$ coherently such that we obtain a path
starting at $-1 \in D_n$ and ending at $1 \in D_n$.

If all curves of $\Delta$ coincide with the corresponding curves of 
$\Gamma$ then the braids $\gamma$ and $\delta$ are equal. 
%(because braids are determined by the isotopy class of their curve diagrams).
Suppose that curves $1,2,\ldots,i-1$ of $\Delta$ agree with
the corresponding curves of $\Gamma$, and the $i$th is the first transverse
one, $1\le i\le n+1$.  This oriented curve has the same 
startpoint as the $i$th curve of $\Gamma$, and first branches off $\Gamma$
either into the upper or the lower component. In the first case we define 
$\delta > \gamma$, in the second $\delta < \gamma$. 
This is well-defined, by proposition \ref{welldef}.\ppar 

{\bf Example}\stdspace All the diagrams in figure \ref{dicurv} are
reduced \wrt each other, and we observe that $1 < \sigma_1 \cdot
\sigma_2\inv < \sigma_1$.
\begin{prop} \sl The relation `$<$' is an ordering, ie if 
$\sigma$, $\gamma$, $\delta$ are braids with $\sigma < \gamma < \delta$ then 
$\sigma < \delta$. \end{prop}
\begin{proof} By the triple reduction lemma \ref{triplem} we can find curve diagrams 
$\Sigma$, $\Gamma$, $\Delta$ of these braids which are all pairwise reduced. 
The statement of the proposition follows immediately: if 
$\Gamma$ branches off $\Sigma$ to the left and $\Delta$ branches off $\Gamma$ 
to the left, then $\Delta$ branches off $\Sigma$ to the left.\end{proof}

\begin{prop} \sl The ordering `$<$' is right invariant.\end{prop}

\begin{proof} Suppose we have two braids $\delta$ and $\gamma$ with $\delta < \gamma$,
and with reduced curve diagrams $\Delta$ and $\Gamma$. 
Let $\sigma$ be a further braid, ie a homeomorphism of $D_n$
fixing $\partial D_n$. We obtain the curve diagrams for $\delta \cdot \sigma$
and $\gamma \cdot \sigma$ by applying $\sigma$ to $\Delta$ and $\Gamma$.
The resulting curve diagrams $\sigma(\Delta)$ and $\sigma(\Gamma)$ are
still reduced, and $\sigma(\Delta)$ still branches off
$\sigma(\Gamma)$ into the lower component of $D_n \basl \sigma(\Gamma)$,
so $\delta \cdot \sigma < \gamma \cdot \sigma$. \end{proof}

Let $\epsilon$ be the trivial braid, with standard curve diagram $E$
(see figure \ref{dicurv}(a)).  We call a braid $\gamma$ {\sl positive}
if $\gamma > \epsilon$, and {\sl negative} if $\gamma < \epsilon$. If
we want to stress that the first $i-1$ curves of $\Gamma$ are parallel
to the corresponding curves of $E$, and the $i$th is the first
non-parallel one, then we say $\gamma$ is {\sl $i$--positive}
respectively {\sl $i$--negative}. Since there is a very similar
concept of $\sigma_i$--positive (see the next section) we shall often
say {\sl geometrically} $i$--positive or negative.  Given two braids
$\gamma$ and $\delta$ such that $\gamma>\delta$ we say $\gamma$ is
(geometrically) {\sl $i$--greater} than $\delta$ if the $i$th curves
are the first non-parallel ones. Any curve diagram in which the first
$i-1$ curves are parallel to the corresponding curves of $E$ is called
$(i-1)$--{\sl neutral}.

We note some simple consequences of right invariance. We have $\gamma
> \epsilon$ if and only if $\epsilon > \gamma \inv$, so the inverse of
a positive braid is negative. If $\gamma > \epsilon$ and $\delta$ is
any braid, then $\gamma\delta > \delta$. ({\it Warning:\/} it need not
be true that $\delta \gamma > \delta$ --- see the next section.)
In particular, the product of positive braids is positive.

% ------------------------------------------------------------------------

\section{Left-consistent canonical form}

In this section we connect our ordering with Dehornoy's
\cite{D.proofnoaxiom}.  The following definition is taken from
\cite{D.proofnoaxiom}. A word of the form
$$
\beta_0 \sigma_i \beta_1 \sigma_i \ldots \sigma_i \beta_k,
$$
where $i\in \{1,\ldots, n-1\}$, and $\beta_0,\ldots \beta_k$ are words
in the letters $\sigma_{i+1}^{\pm1}\ldots \sigma_{n-1}^{\pm1}$ is called a {\sl
$\sigma_i$--positive} word.  A braid is {\sl $\sigma_i$--positive} if it
can be represented by a $\sigma_i$--positive word. A braid is called
{\sl $\sigma_i$--negative} if its inverse is $\sigma_i$--positive.
We shall say that a braid is $\sigma$--positive or negative if it is
$\sigma_i$--positive or negative for some $i$. The following is
the main result from \cite{D.proofnoaxiom}:

\begin{Dthm}
Every braid is precisely one of the following three: $\sigma$--positive, 
or $\sigma$--negative, or trivial.
\end{Dthm}

Dehornoy uses this theorem to define a right-invariant order by 
$\alpha<\beta\iff\alpha\beta\inv$ is $\sigma$--positive.  
We shall prove that this order coincides with the order we 
defined in the last section by showing that the concepts of 
geometrically $i$--positive and $\sigma_i$--positive coincide.  
One way is easy.

\begin{prop}\label{ordagree} A braid which is $\sigma_i$--positive is
geometrically $i$--positive. 
\end{prop}

\begin{proof} A braid which can be represented by a word $\beta \sigma_i$,
where $\beta$ is a word in the letters $\sigma_{i+1}, \ldots,
\sigma_{n-1}$, is geometrically $i$--positive.  To see this think of
the homeomorphism determined by the braid word as a sequence of twists
of adjacent holes around each other: $\beta$ leaves the first $i$
curves untouched and then $\sigma_i$ twists the $i$th hole around the
$(i+1)$st producing a curve diagram in which the $i$th curve moves
into the upper half of $D_n$.  Now by definition, every Dehornoy
positive braid is a product of such words.  The proposition now
follows from the fact that the product of two geometrically $i$--positive
braids is again geometrically $i$--positive.  \end{proof}

The proposition immediately implies part of Dehornoy's theorem:
every braid can take {\it at most} one of the three possible forms.
%since a braid cannot be both geometrically positive and 
%geometrically negative.
To complete the proof that the concepts 
of geometrically $i$--positive and $\sigma_i$--positive coincide
and to recover the remainder of Dehornoy's theorem we shall 
construct a canonical $\sigma_i$--positive form for a given 
geometrically $i$--positive braid.
This is the {\sl left-consistent canonical form} of the braid:

\begin{theorem}[Left-consistent canonical form]\label{lcform}
Let $\gamma$ be a geometrically $i$--positive braid.  Then there is a
canonically defined $\sigma_i$--positive word which represents the same
element of $B_n$.\end{theorem}

We define the {\sl complexity} of a braid $\gamma$ as follows.
Take a curve diagram $\Gamma$ for $\gamma$ which is reduced with
respect to the trivial curve diagram $E$.  Suppose that the first
$j-1$ curves of $\Gamma$ coincide with the first $j-1$ curves of $E$ 
and that the $j$th curve does not.  Let $m\ge0$ be the number of transverse
intersections of $\Gamma$ with $j$th curve of $E$.  The {\sl
complexity of} $\gamma$ is the pair $(j,m)$.  We order complexity
lexicographically with $j$ in reverse order.  Thus $(1,m)$ is more
complex than $(2,n)$ for any $m,n$ whilst $(j,m)$ is more complex than
$(j,n)$ if and only if $m>n$.  The main step in the proof of theorem
\ref{lcform} is the following:  

\begin{prop}\label{mainprop}
Suppose that $\gamma$ is a geometrically $i$--positive braid. 
Then there is a word $\beta$ in the braid
generators $\sigma_{i},\ldots,\sigma_{n-1}$ and their inverses such that 

{\rm(1)}\stdspace $\beta$ contains $\sigma_i\inv$ exactly once

{\rm(2)}\stdspace  $\beta$ does not contain $\sigma_i$

{\rm(3)}\stdspace  $\gamma\beta$ is either geometrically $i$--positive or
geometrically $i$--neutral

{\rm(4)}\stdspace  $\gamma\beta$ has smaller complexity than $\gamma$.

Furthermore there is a canonical choice for $\beta$.

\end{prop}

Theorem \ref{lcform} follows from proposition \ref{mainprop} by 
induction on complexity because, by (4) and induction,
$\gamma':=\gamma\beta$ has a canonical form which by (3) is either
$\sigma_i$--positive or $\sigma_j$--positive or negative for
$j>i$ and then $\gamma'\beta\inv$ is the canonical form for $\gamma$.

\rk{Proof of proposition \ref{mainprop}}
For definiteness we shall deal with the case $i=1$ first.  (We shall
see that the general case is essentially the same as this case.)
So let $\gamma$ be geometrically $1$--positive braid and $\Gamma$
a curve diagram for $\gamma$ which is reduced \wrt the trivial
curve diagram.  We shall define $\beta$ geometrically by sliding one 
particular hole of $D_n$ along a {\sl useful arc}.

Let $E_1\subseteq D_n$ be the $1$st curve of $E$, ie a straight line
 from $-1$ to the leftmost hole of $D_n$, excluding this hole.  We
define a {\sl useful arc} to be a segment $b$ of some curve of
$\Gamma$ starting at some point of $E_1$ (possibly $-1$), and ending
at some hole of $D_n$ other than the leftmost one such that
 
\begin{itemize}
\item[-] the interior of $b$ does not intersect $E_1$,
\item[-] an initial segment of the arc $b$ lies in the upper half of the disk, 
ie the intersection of a neighbourhood of $E_1$ with the interior of $b$ consists of a 
line segment in the upper component of $D_n\basl E$.
\end{itemize}

\begin{figure}[htb] 
\cl{ \relabelbox \epsfbox{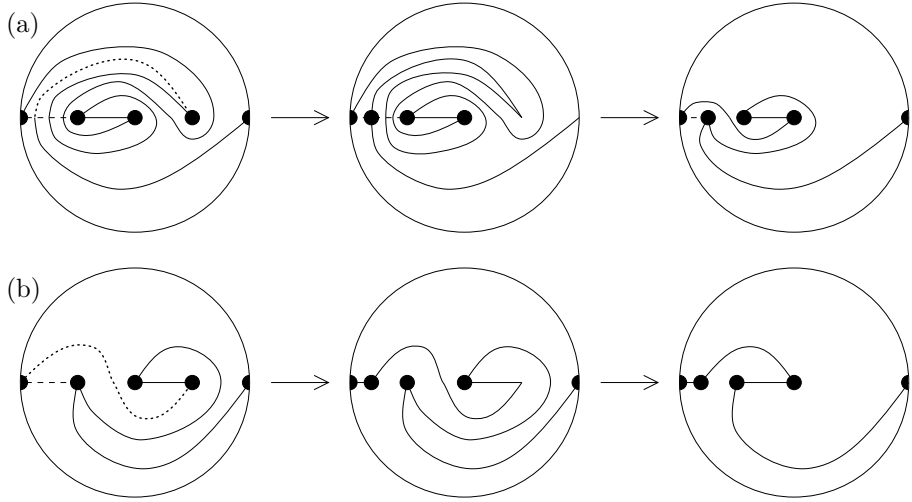} \relabel{(a)}{\small (a)}
\relabel{(b)}{\small (b)} \endrelabelbox }
% exported at 60pc portrait
\caption{Slide of a hole along a useful arc, followed by a reduction} 
\label{exrdisks}
\end{figure}

Suppose that $\Gamma$ contains useful arcs. Then each of them has
precisely one point of intersection with $E_1$ and we call the one
whose intersection point is leftmost the {\sl leftmost useful arc}.
Let $b$ be the leftmost useful arc.  If $b$ starts in the interior of
$E_1$, then we can slide the hole of $D_n$ at the endpoint of $b$
along $b$ and back into $E_1$.  If $b$ starts at $-1$, then we push a
small initial segment of $b$ into $E_1$, and then perform the slide of
the hole of $D_n$ (see figure \ref{exrdisks} where $b$ is dotted). In
either case we obtain a curve diagram $\Gamma'$ representing a braid
$\gamma'$. Now $\Gamma'$ need not be reduced \wrt $E$.  But notice
that $\gamma'$ has lower complexity than $\gamma$ since the new $E_1$
now stops at the intersection of $b$ with the old $E_1$ and hence
there are fewer intersections with $\Gamma'$ even before reduction.

The movement of the hole of $D_n$ along $b$ defines a braid $\beta$ 
on $n$ strings, with $\gamma' = \gamma\beta$.  Furthermore we can 
decompose $\beta$ as a canonical word in the generators $\sigma_j$
by writing down the appropriate $\sigma_j$ or $\sigma_j\inv$ 
whenever the hole passes over or under another hole.  But, by definition
of useful arc, the hole only passes once over or under the leftmost 
hole  and it passes over and to the left and hence the word that we 
read contains $\sigma_1\inv$ only once and does not contain $\sigma_1$.

Therefore to prove case $i=1$ of proposition \ref{mainprop} 
it remains to prove the following two claims:

{\bf Claim 1}\stdspace The diagram contains a useful arc.

{\bf Claim 2}\stdspace The diagram $\Gamma'$ obtained by sliding a 
hole of $D_n$ along the leftmost useful arc is either $1$--positive 
or $1$--neutral, but not $1$--negative.
 
\begin{figure}[hbt] 
\cl{
\relabelbox
\epsfbox{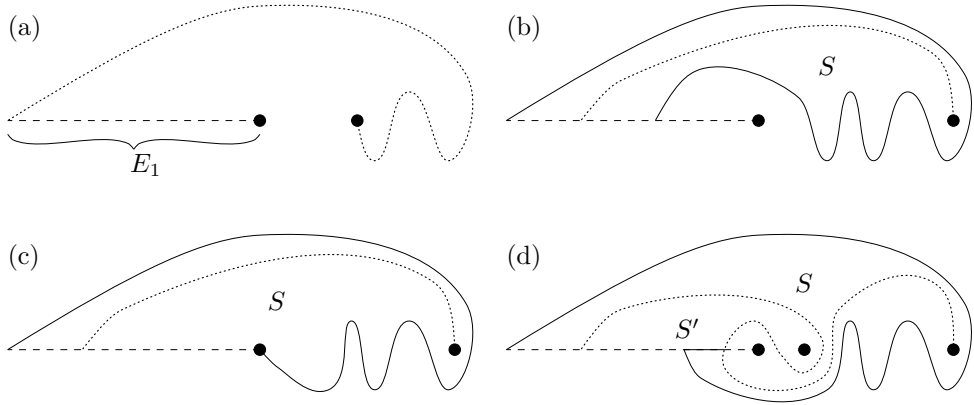}
\relabel{(a)}{\small (a)}
\relabel{(b)}{\small (b)}
\relabel{(c)}{\small (c)}
\relabel{(d)}{\small (d)}
\relabel{S}{\small $S$}
\relabel{S+}{\small $S$}
\relabel{S++}{\small $S$}
\relabel{S'}{\small $S'$}
\relabel{E0}{\small $E_1$}
\endrelabelbox
}
% exported at 48pc as portrait
\caption{How to find a useful arc} \label{excerd}
\end{figure}

To prove claim 1, we consider the first curve of $\Gamma$ starting at
$-1$.  If it ends in a hole other than the leftmost one and does not
intersect $E_1$ then it is a useful arc (figure \ref{excerd}(a)).
Otherwise we consider the closed curve in $D^2$ starting at $-1$,
along the first curve of $\Gamma$, up to its first intersection with
the closure of $E_1$ in $D^2$, and then back in a straight line to the
point $-1$. This curve bounds a disk $S$ in $D^2$, which may be of
three different types: $\Gamma$ hits $E_1$ either from above, or from
below, or in the leftmost hole of $D_n$ (see figure
\ref{excerd}(b),(d),(c)). In cases (b) and (c) we note that since
$\Gamma$ and $E$ are reduced, at least one hole of $D_n$ must lie in
the interior of $S$. Moreover, all holes of $D_n$ are connected by
curves of $\Gamma$, so there exists a curve of $\Gamma$ connecting one
of the holes in $S$ to one of the holes outside $S$ or the point $1
\in D_n$. The first component of the intersection of this curve with
$S$ is a useful arc.

In case (d) we walk along the oriented curve in $D^2$ starting at $-1$,
along the curves of $\Gamma$. We write down the symbol $+$ whenever we 
hit $E_1$ from below (or at $-1$), and $-$ if we hit $E_1$ from above 
or in the leftmost 
hole of $D_n$. The sequence starts with a $+$, and since the curve has to leave
the disk $S$ it must contain a $-$. It follows that the string $+-$ must
occur in the sequence; it represents an arc which, together with a segment 
of $E_1$, bounds a disk $S'$ in $D^2$. See figure \ref{excerd}(d): 
$S'$ is bounded
by part of the dotted arc between two intersections with $E_1$ and
part of $E_1$.  Since $\Gamma$ and $E$ are reduced, 
$S'$ contains a hole other than the leftmost one in its boundary or in its 
interior. In the first case, a segment of top (dotted) boundary of $S'$
is a useful arc; 
in the second case the disk $S'$ is of the type indicated in figure 
\ref{excerd}(b) or (c), so there is a useful arc inside $S'$.
This finishes the proof of claim 1.
\ppar

To prove claim 2, we distinguish two cases: either the leftmost useful
arc $b$ starts at the point $-1$, or it starts at some point in the
interior of $E_1$.  In the first case (eg figure \ref{exrdisks}(b)),
the curve diagram $\Gamma'$ obtained by sliding a hole along $b$ to
near $-1$ is $1$--neutral.

\begin{figure}[hbt] 
\cl{
\relabelbox
\epsfbox{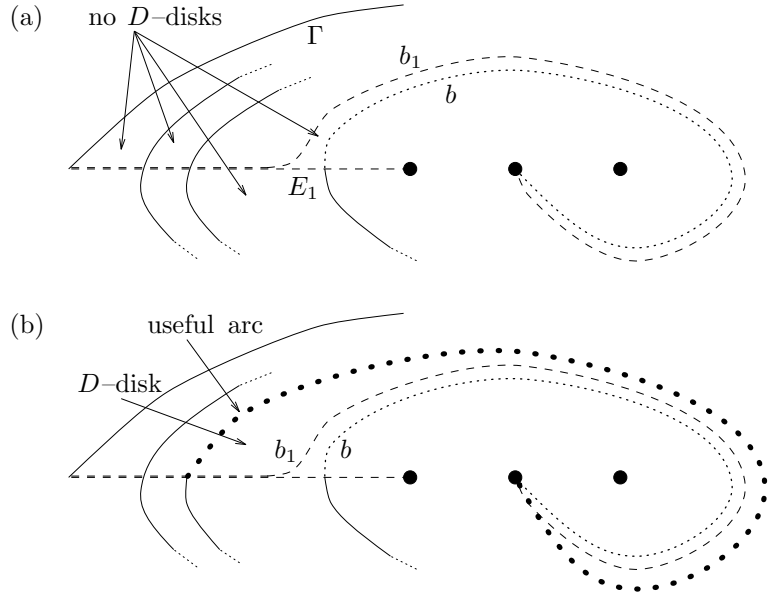}  %%%%  comma , deleted for Dos reasons
\relabel{(a)}{{\small (a)}}
\relabel{(b)}{{\small (b)}}
\relabel{no}{{\small no $D$--disks}}
\relabel{D}{}
\relabel{-disks}{}
\relabel{Ddisk}{{\small $D$--disk}}
\relabel{E1}{\small $E_1$}
\relabel{b}{\small $b$}
\relabel{b0}{\small $b_1$}
\relabel{b'}{\small $b$}
\relabel{b0'}{\small $b_1$}
\relabel{useful arc}{{\small useful\ \,arc}}
\relabel{G}{{\small $\Gamma$}}
\endrelabelbox
}
\vspace{-1mm}
% exported at 55pc as portrait
\caption{There are no $D$--disks between $b_1$ and $\Gamma$}
\label{bb_0const}
\end{figure}

In the second case (figure \ref{exrdisks}(a)) the curve diagram
$\Gamma'$ is $1$--positive, as we now prove. We recall that we had
$\gamma' = \gamma \beta$, where $\beta$ represents the slide of a hole
along the leftmost useful arc $b$. We observe that we can construct a
curve diagram of the braid $\beta \inv$ such that the first curve
$b_1$ of the diagram is a line segment in $E_1$ from $-1$ almost all
the way to $E_1 \cap b$, followed by an arc parallel and close to the
arc $b$, and finally running into the same hole as $b$. The
construction of the arc $b_1$ is illustrated in figure
\ref{bb_0const}(a).

Next we examine the possible reductions of $\Gamma$ \wrt this arc $b_1$. 
If there was a $D$--disk of type (b) whose boundary contained the arc $b$,
(ie to the right of $b_1$ in figure \ref{bb_0const})  
then cutting off the strip bounded by $b$, $b_1$ and $E_1$ would yield a 
$D$--disk of type (a) of $\Gamma$ \wrt $E$ (see figure \ref{bb_0const}(a)). 
This is impossible by hypothesis.
If there was a $D$--disk of type (b) whose boundary contained a final segment
of the arc $b_1$ and a segment other than $b$ of a curve of $\Gamma$, 
(ie to the left of $b_1$ in figure \ref{bb_0const}) then 
this segment would be a useful arc intersecting $E_1$ more to the left
than $b$ (figure \ref{bb_0const}(b)), which is also impossible. 
Finally, any
%, so no such $D$--disk can exist. 
$D$--disk of type (a) of $\Gamma$ \wrt $b_1$ would also be a $D$--disk
of $\Gamma$ \wrt $E_1$. So there are no $D$--disks between $b_1$ and $\Gamma$.
By remark \ref{onecurv} it follows that we can reduce
the curve diagram of $\beta\inv$ \wrt $\Gamma$ without touching its first
curve $b_1$. We can now observe that $\gamma$ is $1$--greater than 
$\beta\inv$, ie $\gamma'$ is $1$--positive, as claimed. This completes the 
proof of claim 2. 

Finally we turn to the case when $i$ may not be 1.  In this case 
the first $i-1$ holes are lined up near $-1$ on the real axis.
The same argument as in the case $i=1$, only with the $i-1$st hole
and the line segment $E_i$ playing the role previously played by
$-1$ and $E_1$ respectively, completes the proof of the general case.
%we slide the first $i-1$ holes along the real axis to near $-1$ and 
%renumber so that the $i$th hole is now hole 1.  We then use
%the argument of the case $i=1$ and note that the end of $E_1$ near 
%$-1$ is not disturbed, so restoring the correct numbers to the holes
%we have proved the general case.
\qed\ppar

The proof of theorem \ref{lcform} provides an explicit algorithm for 
converting a braid into its left-consistent canonical form.  In the 
appendix we give a formal version of this algorithm using cutting sequences.

\rk{Remark}
The order on the braid group has the property that inserting
a generator $\sigma_i$ anywhere in a braid word makes the braid larger.
A proof of this fact, in the spirit of this paper, is given in 
\cite{bert}.  This property is equivalent to the statement that the 
order extends the {\sl subword order} defined by Elrifai and Morton 
\cite{elfir} and an algebraic proof has been given by Laver \cite{laver}.

% -----------------------------------------------------------------

\section{Counterexamples}

We shall call a braid word {\sl $\sigma$-consistent} (Dehornoy in 
\cite{D.algorithm} calls it {\sl reduced}) if it is $\sigma$-positive,
$\sigma$-negative or trivial. We have seen in the previous chapter that 
every braid has at least one $\sigma$-consistent representative. The aim 
of this chapter is to disprove some plausible-sounding but overoptimistic 
conjectures about the ordering and about $\sigma$-consistent representatives 
of braids. 

%blabla
%{\small {\bf Abstract } We prove that every element of the braid groups
%$B_2$, $B_3$ has a simultaneously shortest and left consistent representative,
%whereas the analogue statement for $B_n$, $n\geqslant 4$, is false. We also
%give an example that the left consistent canonical form of a braid does not
%necessarily have the minimal number of occurences of the main generators 
%among all left consistent representatives. } rubbish...
%\ppar

%
\begin{figure}[hbt] 
\cl{
\relabelbox
\epsfbox{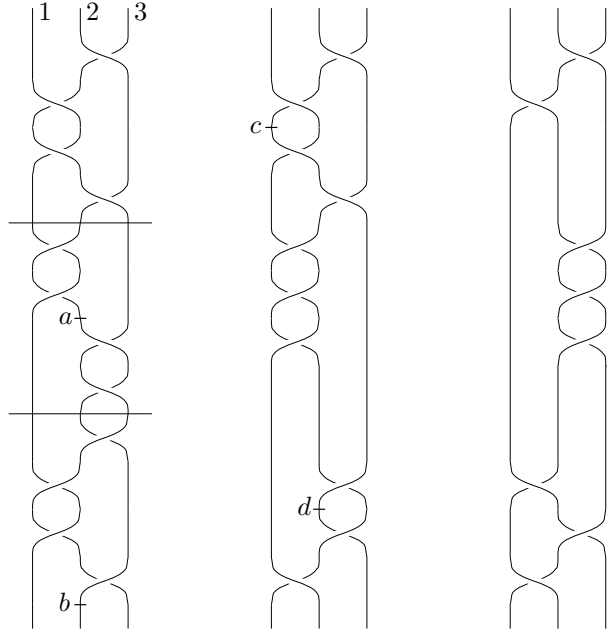} 
\relabel{a}{\small $a$}
\relabel{b}{\small $b$}
\relabel{c}{\small $c$}
\relabel{d}{\small $d$}
\relabel{1}{\small 1}
\relabel{2}{\small 2}
\relabel{3}{\small 3}
\endrelabelbox
}
% exported at 50pc as portrait
\caption{Equivalence of a conjugate of a positive pure braid with a visibly 
negative braid} \label{pureconj}
\end{figure}
\sh{Left invariance on the pure braid group}

Because the {\it pure} braid group has an
ordering which is simultaneously left and right invariant
\cite{lrinvariant}, it would be tempting to think that the geometric 
ordering is left and right invariant when restricted to the pure braid
group. 
However this is equivalent to saying that a pure positive
braid, when conjugated by any pure braid, is again positive and the
example in figure \ref{pureconj} shows this to be false.
In $B_3$, the braid group on three strings, we conjugate the pure positive 
braid $\sigma_1^2 \sigma_2^{-2}$ by the pure braid $\sigma_2\
\sigma_1^2\ \sigma_2$. The figure shows the equivalence of the 
resulting braid with the $\sigma$--negative braid 
$\sigma_2\inv \sigma_1\inv \sigma_2^3\ \sigma_1\inv \sigma_2\ \sigma_1\inv$.
%$\sigma_1\inv \sigma_2 \sigma_1\inv \sigma_2^3\ \sigma_1\inv \sigma_2\inv$.
We are moving first the string segment $\overline{ab}$ and then the segment 
$\overline{cd}$ `over' the braid `to the left of the braid'.

\sh{Simultaneously shortest and $\sigma$-consistent representatives}

For any element $b$ of the braid group $B_n$ ($n\geqslant 2$),
%$$ B_n=\langle \sigma_1,\ldots, \sigma_{n-1}\ | \ \sigma_i\sigma_j=
%\sigma_j\sigma_i \hbox{ \ if \ }|i-j| \geqslant 2,\ \sigma_{i+1}\sigma_i
%\sigma_{i+1}=\sigma_i\sigma_{i+1}\sigma_i \rangle,$$ where $n\geqslant 2$, 
there are two ways to represent $b$ by a particularly simple word $w$ in 
the letters $\sigma_1^{\pm 1}$,$\ldots$,$\sigma_{n-1}^{\pm 1}$.

\begin{figure}[htb] 
\cl{
\relabelbox
\epsfbox{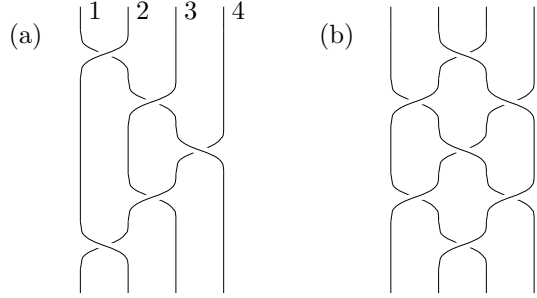} % exported at 50pc as portrait
\relabel{1}{\small $1$}
\relabel{2}{\small $2$}
\relabel{3}{\small $3$}
\relabel{4}{\small $4$}
\relabel{(a)}{\small (a)}
\relabel{(b)}{\small (b)}
\endrelabelbox
}
\caption{The equivalent braids $\sigma_1 \sigma_2 \sigma_3\inv \sigma_2 
\sigma_1\inv$ and $\sigma_2\inv \sigma_3\inv \sigma_1 \sigma_2\inv \sigma_1
\sigma_3 \sigma_2$}  \label{twobr}
\end{figure}

(1) $b$ can be represented by a word which is as short as possible.
For instance, we shall see later that the word $w_1=\sigma_1 \sigma_2 
\sigma_3\inv \sigma_2 \sigma_1\inv$ is a shortest possible representative
of a braid in $B_4$ (see figure \ref{twobr}(a)).

(2) $b$ can be represented by a $\sigma$-consistent word. 
%(Dehornoy in \cite{D.algorithm} calls it a {\sl reduced} word). That means,
%$w$ is either the nullstring, or the {\sl main generator}, ie the generator 
%with lowest index occuring in $w$, occurs only positively or only negatively. 
%(The fact that every braid posesses such a left consistent representative was 
%first proved in \cite{D.proofnoaxiom}, a geometric proof is given in 
%\cite{bror1}.) 
For instance, in the braid word $w_2=\sigma_2\inv \sigma_3\inv \sigma_1 
\sigma_2\inv \sigma_1 \sigma_3 \sigma_2$, which represents the same element 
of $B_4$ as $w_1$, the letter $\sigma_1$ occurs only with positive exponent,
see figure \ref{twobr}(b).

\begin{theorem} Every element of $B_n$ for $n=2,3$ has a simultaneously
shortest and $\sigma$-consistent representative. By contrast, there are
braids in $B_n$ for $n\geqslant 4$ all of whose $\sigma$-consistent 
representatives have non-minimal length. \end{theorem}
\begin{proof} The case $n=2$ is obvious. The case $n=3$ follows from the fact that
in $B_3$ Dehornoy's handle-reduction algorithm \cite{D.algorithm} never
increases the length of a braid word, and hence turns any 
shortest representative of a given braid into a simultaneously shortest and
$\sigma$-consistent one. 

For the case $n\geqslant 4$ it suffices to prove that the braid 
$b:=\sigma_1 \sigma_2 \sigma_3\inv \sigma_2 \sigma_1\inv \in B_4$ 
%illustrated in 
(figure \ref{twobr}) has length $5$, while every $\sigma$-consistent 
representative has more than five letters.

To see that every representative has at least five letters we note that
the image of $b$ under the natural homomorphism $B_n \to S_n$, from the braid
group into the symmetric group, is the permutation $(14)$. This permutation
cannot be written as a product of less than five adjacent transpositions. 
The result follows.

We now assume, for a contradiction, that there exists a five-letter 
representative which is also $\sigma$-consistent. This would be a braid on 
four strands with the following properties:
\begin{itemize}
\item[(i)] its image under the natural map $B_4 \to S_4$ is $(14)$,
\item[(ii)] it has five crossings (ie it is a word with five letters),
\item[(iii)] if we denote by $c(i,j)$ ($i,j \in \{1,\ldots, 4\}$) the
algebraic crossing number of the $i$th and the $j$th string, then the braid 
must satisfy 
$c(1,2)=1$, $c(1,3)=1$, $c(1,4)=-1$, $c(2,3)=0$, $c(2,4)=-1$, $c(3,4)=1$,
\item[(iv)] it may contain the letter $\sigma_1$, but not $\sigma_1\inv$
(note that there exists a representative of $b$ in which $\sigma_1$ occurs 
only positively, so there can't exist a consistently negative one).
%, see \cite{Larue,D.proofnoaxiom,bror1}).
\end{itemize}
\begin{figure}[htb] 
\cl{
\relabelbox
\epsfbox{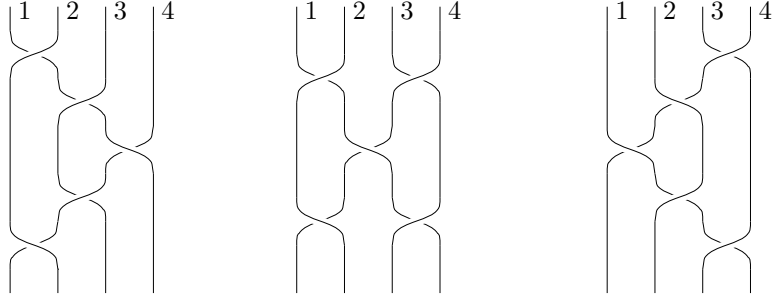} % exported at 50pc as portrait
\relabel{1}{\small $1$}
\relabel{2}{\small $2$}
\relabel{3}{\small $3$}
\relabel{4}{\small $4$}
\relabel{1'}{\small $1$}
\relabel{2'}{\small $2$}
\relabel{3'}{\small $3$}
\relabel{4'}{\small $4$}
\relabel{1''}{\small $1$}
\relabel{2''}{\small $2$}
\relabel{3''}{\small $3$}
\relabel{4''}{\small $4$}
\endrelabelbox
}
\caption{Three candidates for short $\sigma$-consistent representatives of $b$}
\label{threepos}
\end{figure}
There are only three braids satisfying (i) - (iii), pictured in figure 
\ref{threepos}, and we observe that none of them satisfies (iv). It follows 
that no $\sigma$-consistent representative of $b$  with only five crossings 
exists.
\end{proof}

\sh{Minimal number of occurrences of the main generator}

We define the {\sl main generator} of a braid word to be the generator
with lowest index occurring in the word. 
It is tempting to think that sliding holes along {\it leftmost} useful arcs,
as in the left consistent canonical form, is the most efficient way of 
reducing the number of intersections between the curve diagram and the
line segment $E_1$. This, however, is wrong:
\begin{theorem} There are braids whose left consistent canonical form
does not have the minimal number of occurrences of the main generator 
among all $\sigma$-consistent representatives. \end{theorem}
%

%We finish by stating the following related conjecture, which is due to 
%C Rourke and M T Greene:
%The following related, very plausible, conjecture was made independently
%by M.~T.~Greene and C.~Rourke \cite{perscom}:
%
%{\bf Conjecture } {\sl Among all left consistent words representing a given 
%braid, the left consistent canonical form (as defined in \cite{bror1}) has 
%the minimal number of occurrences of the main generator.}
%
%{\bf Theorem 2 } {\sl The conjecture is false - the braid $\Delta^3 \in B_3$,
%where $\Delta=\sigma_2\sigma_1\sigma_2$ is a counterexample.}

\begin{figure}[htb] 
\cl{
\relabelbox
\epsfbox{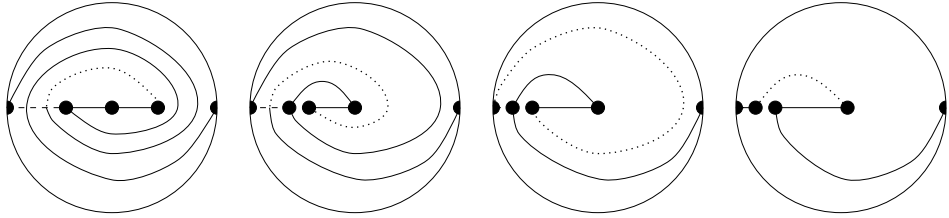} % exported at 55pc as portrait
%\relabel{1}{\small $1$}
\endrelabelbox
}
\caption{The left consistent canonical form of $\Delta^3$ is 
$\sigma_2\sigma_1\sigma_2\sigma_2\sigma_1\sigma_2\sigma_2\sigma_1\sigma_2$} 
\label{delta3cd}
\end{figure}

\begin{proof} We shall show that the braid $\Delta^3$, where 
$\Delta=\sigma_2\sigma_1\sigma_2$, has this property.
Note that $\Delta$ is just a half-twist, so $\Delta^2$ generates the 
commutator subgroup of $B_3$.

We have $(\sigma_2\sigma_1\sigma_2)^3=\sigma_2\sigma_2\sigma_1
\sigma_2\sigma_2\sigma_2\sigma_1\sigma_2\sigma_2$, so the braid can be
represented by a $\sigma$-consistent word in which the main generator 
$\sigma_1$ occurs only twice. However, as is easy to check with the help 
of figure \ref{delta3cd}, the left consistent canonical form of the braid 
is the word $(\sigma_2\sigma_1\sigma_2)^3$, which contains the main 
generator $\sigma_1$ three times. \end{proof}

\sh{Local indicability}

We are indebted to Stephen P Humphries and Jim Howie for pointing out
the following. A group is called {\sl locally indicable} if every
finitely generated subgroup has a nontrivial homomorphism to the
integers. It was proved by Burns and Hale \cite{BuHa} that locally
indicable groups are right-orderable, but it took almost two decades
until G Bergman \cite{Bergman} found an example of a group which is
right-orderable but not locally indicable; ie the class of locally
indicable groups is {\it strictly} contained in the class of
right-orderable groups. We can now give further examples:

\begin{theorem} The braid group $B_n$ for $n\geqslant 5$ is right orderable
but not locally indicable. \end{theorem}

\begin{proof}
It remains to show that $B_n$ is not locally indicable. The commutator 
subgroup $B_n'$ of $B_n$ is finitely generated, and for $n\geqslant 5$ the 
first and second commutator subgroups coincide: $B_n'=B_n''$ 
(see \cite{GoLin}). It follows that the abelianization of $B_n'$ is trivial, 
so $B_n'\subset B_n$ has no nontrivial homomorphism to $\Z$.
\end{proof}

% --------------------------------------------------------------------------

\section{Automatic ordering}

Define a right-invariant ordering to be {\sl automatic} if it can be
determined by a finite-state automaton.  In this section we shall see
that the ordering on the braid group is automatic.

This is proved by comparing the order on the braid group as defined
in section 3 with Mosher's automatic structure \cite{Mosher1, Mosher2}.  
This comparison gives more.
Define a group to be {\sl order automatic} if it is both automatic
\cite{Ep} and right-orderable and such that there is a finite
state automaton which detects the order from the automatic normal
forms.  To be precise, there exists an automatic structure and a
finite state automaton, which, given two normal forms for the
automatic structure, will decide which represents the greater group
element.

\begin{theorem}\label{autoorder}
The braid group $B_n$ is order automatic.\end{theorem} 

\begin{remark}\label{strongalg} The algorithm to decide which of two 
given normal forms is the greater takes linear time in the length 
of the normal form.  Using results from Epstein et al \cite{Ep} we 
deduce: \end{remark}

\begin{cor}\label{quadtime}There is a quadratic-time algorithm to 
decide which of two elements of $B_n$ (presented in terms of standard 
braid generators) is the greater.\end{cor}

Full details of the proof of these results can be found in
\cite{Rourke-Wiest}.  Here we shall give a short proof of theorem
\ref{autoorder} which yields only a quadratic time algorithm to 
order normal forms which is nevertheless sufficient to imply corollary
\ref{quadtime}. 

In \cite{Mosher1,Mosher2} Mosher constructs normal forms for elements
of mapping class groups by {\sl combing} triangulations (and hence
proves that mapping class groups are automatic). We shall need to
sketch Mosher's normal form in the special case of the braid group.

\begin{figure}[hbt] 
\cl{
\relabelbox\small
\epsfxsize 2truein
\epsfbox{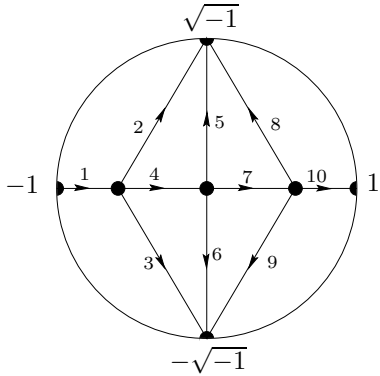}
\def\ss#1{$\scriptstyle#1$}
\adjustrelabel <0pt, 1pt> {1}{\ss1}
\relabel{2}{\ss2}
\relabel{3}{\ss3}
\adjustrelabel <1pt, 1pt> {4}{\ss4}
\relabel{5}{\ss5}
\relabel{6}{\ss6}
\adjustrelabel <-2pt, 1pt> {7}{\ss7}
\relabel{8}{\ss8}
\relabel{9}{\ss9}
\adjustrelabel <1pt, 0pt> {10}{\ss{10}}
\relabel{p1}{1}
\relabel{m1}{$-1$}
\relabel{i}{$\sqrt{-1}$}
\adjustrelabel <-10pt, -2pt> {mi}{$-\sqrt{-1}$}
\endrelabelbox
}
\caption{The base triangulation} \label{basetr}
\end{figure}

We define the {\sl base triangulation} $B$ of $D_n$ to have vertices
at the $n$ missing points and at the four {\sl boundary vertices},
$\pm1$ and $\pm\sqrt{-1}$.  The edges of $B$ comprise the four arcs of
$\pa D_n$ joining pairs of boundary vertices, $n+1$ edges along the
real axis and $2n$ edges joining $\pm\sqrt{-1}$ to the real vertices
not $\pm1$, see figure \ref{basetr}.  We order and orient the edges as
indicated.

An {\sl allowable triangulation} of $D_n$ is a triangulation with the
same vertex set.  We identify two allowable triangulations if they
differ by a vertex fixing isotopy.  A {\sl triangulation class} is a
set of boundary fixing isomorphism classes of allowable
triangulations.  Ie two triangulations are in the same class if they
differ by an element of the braid group.

We now consider the groupoid $\G$ which has for objects the set of
triangulation classes of $D_n$ and for morphisms the set of ordered
pairs $(T,T')$ of allowable triangulations, where $(T,T')$ is
identified with $(h(T),h(T'))$ if $h\in B_n$.  The morphism goes from
the class of $T$ to the class of $T'$.  If $T$ and $T'$ are in the
same class, then there is a unique boundary fixing isomorphism from
$T'$ to $T$ up to isotopy, ie an element of $B_n$.  This determines an
isomorphism between the vertex group of $\G$ and the braid group
$B_n$.  (Note that for this isomorphism, and for compatibility with
Mosher's conventions, we need to replace the {\sl algebraic}
convention for multiplication in the braid group, described in section
1, by the opposite {\sl functional} convention, ie
$\Phi\Psi:=\Phi\circ\Psi$.  The functional convention is used
throughout this section; the algebraic convention is used in all
other sections and in the appendix.)

\sh{Combing}

We consider a particular type of morphism in $\G$.

\begin{figure}[hbt] 
\cl{
\relabelbox\small
\epsfxsize 2truein
\epsfbox{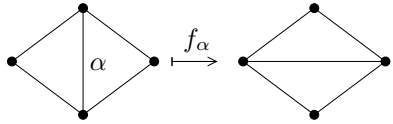}
\relabel{alpha}{$\alpha$}
\adjustrelabel <0pt, 4pt> {f_a}{$f_\alpha$}
\endrelabelbox
}
\caption{Flipping an edge} \label{flipedg}
\end{figure}

\rk{Definition}{\sl Flipping an edge}\stdspace  An edge $\alpha$ adjacent
to two triangles $\delta$ and $\delta'$ is removed (to form a square
of which $\alpha$ is a diagonal) and then the square is cut back into
two triangles by inserting the opposite diagonal.  We call this
morphism ``flipping $\alpha$'' and denote it $f_\alpha$, see figure
\ref{flipedg}.

Every morphism $q=(B,T)$ in $\G$ from the base vertex to another
vertex is a product of a canonical sequence of flips.  To see this,
picture $q$ as given by superimposing $B$ and $T$, and {\sl comb} $T$
along $B$.  To be precise, first reduce $T$ \wrt $B$ and then consider
edge $1$ of $B$.  Suppose that, starting at $-1$, edge one crosses
edge $\alpha$ of $T$. Flip $\alpha$.  Repeat until there are no more
crossings of edge 1 with $T$. (The fact that this process is finite
follows from a simple counting argument: one counts the number of 
intersections of edge 1 with $T$, except with the next edge of $T$ which
is to be flipped.  For more detail here see \cite[pages 321--322]{Mosher2}.)  
Now do the same for edge 2 starting
at the non-boundary vertex and continue in this way, using the
ordering and orientation of edges of $B$ indicated in figure
\ref{basetr}, until $T$ has been converted into a copy of $B$.

The {\sl Mosher normal form} of $q$ is the inverse of the sequence of
flips described above.%
\footnote{Strictly speaking the Mosher normal form is not this
flip sequence, which only defines an asynchronous automatic structure,
but is derived from it by clumping flips together into blocks called
``Dehn twists'', ``partial Dehn twists'' and ``dead ends'' (see
\cite{Mosher2} pages 342 et seq).  This technicality does not affect
any of the results proved here or in \cite{Rourke-Wiest}.  We prove
that order can be detected in linear time from the flip sequence.
Since the clumped flip sequence can be unclumped in linear time, this
implies that order can be detected in linear time from the strict
Mosher normal form.}
Notice that unlike the general case described in \cite{Mosher2}, $q$
is completely characterised by the sequence of flips, there is no need
to carry the labelling of $T$ along. %in the combing sequence.  In
particular, there is no relabelling morphism required here.  (This is
because $\pa D_n$ is fixed throughout.)

\sh{Detecting order from the Mosher normal form}

To see the connection with order, consider an element $(B,T)$ of the
vertex group at the class of $B$.  There is an element $g\in B_n$ (a
homeomorphism of $D_n$ fixing $\pa D_n$) unique up to isotopy
carrying $T$ to $B$.  Conversely given $g\in B_n$ the corresponding
triangulation pair is $(B,g\inv B)$.

We observe that if we comb $B$ along $g(B)$ this is combinatorially
identical to combing $g\inv(B)$ along $B$. We call the sequence of
flips defined by this combing the {\sl combing sequence of $g$}.
(The reverse of the combing sequence is the Mosher normal form of $g$.)

The curve diagram of $g$ is part of the triangulation $g(B)$ namely
the edges numbered $1,4,7,\ldots,3n+1$.  Suppose that $g$ is
$i$--positive, then $g$ can be assumed to fix the first $i-1$
of these edges (ie $1,4, \ldots,3i-5$) and then, after reduction, 
can be assumed to fix the corresponding outlying edges (ie 
$2,5,\ldots,3i-4$ and $3,6,\ldots, 3i-3$).  But edge $3i-2$ is carried
into the upper half of $D_n$ and must meet edge $3i-1$ of $B$.  Thus
the first flip in the combing sequence of $g$ is $f_{3i-1}$,
ie flip the edge numbered $3i-1$.  Similarly if $g$ is $i$--negative
then the first flip in the combing sequence is $f_{3i}$.  We have
proved the following:

\begin{algo}\label{quickalg}
\rm (To decide from the Mosher normal form
whether a braid element is $i$--positive or negative and provide 
the correct value of $i$)\sl\nl
Inspect the combing sequence (the reverse of the Mosher normal form).
The first flip is either $f_{3i-1}$ for some $i$ or $f_{3i}$ for
some $i$.  In the first case the braid is $i$--positive and in the
second it is $i$--negative.\end{algo}

This algorithm is visibly executable by a finite-state automaton and
linear in the length of the normal form of $g$.  Theorem \ref{autoorder}
and corollary \ref{quadtime} follow from general principles.  To
decide the relative order of two elements $\alpha$ and $\beta$ we 
compute the normal form of $\alpha\beta\inv$ --- this can be done by
a finite-state automaton and takes quadratic time, see \cite{Ep} ---
and then apply algorithm \ref{quickalg}.

\rk{Final remarks}
(1)\stdspace 
We have proved that there is a quadratic time algorithm to decide 
the relative order of two braid words.  In \cite{D.algorithm}
Dehornoy presents an algorithm which does this in practice and
is apparently extremely fast --- however his formal proof that this
algorithm works only provides an exponential bound on time.  The
algorithm presented here is implementable since the whole Mosher
program can be implemented, see \cite{Mosher.practical}.  Note that
in the  appendix we present another algorithm based on cutting sequences.

(2)\stdspace There is a far stronger connection between the Mosher
normal form and the order on $B_n$ than presented here.  The relative
order of two elements can be detected from their combing sequences by
inspecting just the first four differences in the sequences (and this
proves remark \ref{strongalg}).  Full details here are to be found in
\cite{Rourke-Wiest}.

% --------------------------------------------------------------------------

\sh{Addresses:}

{\small
R.\ Fenn:\stdspace {\it School of Mathematical Sciences, 
University of Sussex, Falmer, Brighton BN1 9QH, UK}\stdspace 
{\tt R.A.Fenn@sussex.ac.uk}

M.\ T.\ Greene:\stdspace {\it Radan Computational, Ensleigh House, 
Granville Road, Bath BA1 9BE, UK}\stdspace {\tt Michael.Greene@uk.radan.com}

D.\ Rolfsen:\stdspace {\it Department of Mathematics,
University of British Columbia, Vancouver, B.C.\ Canada V6T 1Z2}
\stdspace {\tt rolfsen@math.ubc.ca}

C.\ Rourke:\stdspace{\it Mathematics Institute, University of Warwick, 
Coventry CV4 7AL, UK} \stdspace {\tt cpr@maths.warwick.ac.uk}

B.\ Wiest:\stdspace{\it CMI, Universit\'e de Provence, 
13453 Marseille cedex 13, France,}\newline {\tt bertw@gyptis.univ-mrs.fr} 
}

% ------------------------------------------------------------------------

\appendix\small
\section{Appendix: Cutting sequences}  \label{wordpr}
\setlength{\parskip}{5pt plus3pt minus3pt}

In this appendix we define a unique {\sl reduced cutting sequence} for
a braid.  We give implementable algorithms to read the reduced cutting 
sequence from the braid word, to decide order from the cutting 
sequence and to put a braid, given in terms of standard twist generators, 
into its left-consistent canonical form. 

\sh{Cutting sequences and curve diagrams}

A {\sl cutting sequence} is a finite word $\chi$ in the letters
$0,\ldots, n$, $\underline{0},\ldots$, $\underline{n+1}$, $\up$ and $\down$ 
such that 

\begin{itemize}
\item[(i)] $\chi$ starts with $\underline{0}$ and ends with $\underline{n+1}$,
\item[(ii)] each of the letters $\underline{0},\ldots, \underline{n+1}$
occurs precisely once in $\chi$,
\item[(iii)] in the word $\chi$ numbers and arrows alternate, with the single
possible exception that strings of the form $\underline{i}\ \underline{i+1}$
or $\underline{i+1}\ \underline{i}$ ($i=0,\ldots,n$) may occur.
\end{itemize}

Consider now a curve diagram $\Gamma$. It consists of three types of
subcurves: curves in the upper half plane, curves in the lower half
plane, and straight line segments in the real line. Note that curves 
in the upper or lower half plane may be replaced by semicircles since
they are determined by their end points. For convenience we rescale the 
curve diagram so that it goes from $0$ to $n+1$ and the $n$ holes are 
the integers $1,2,\ldots,n$.

Going along $\Gamma$ we can read off a cutting sequence, by reading an
$\up$ or $\down$ for every curve in the upper or lower half plane respectively,
an $\underline{i}$ ($i\in \{0,\ldots, n+1\}$) for every intersection with
the integer $i$ in the real line (so underlined integers correspond
to holes), and an $i$ for every intersection with the real interval
$(i, i+1)$. It is easy to check that a word obtained in this way is indeed
a cutting sequence.

For example the curve diagram representing $\sigma_1$ in figure \ref{dicurv}
is coded as $\underline{0} \up \underline{2}\ \underline{1} \down \underline{3}
\ \underline{4}$, whereas $\sigma_1 \sigma_2\inv$ is coded $\underline{0} \up
1 \down \underline{3} \down \underline{1} \down 3 \up \underline{2} \up
\underline{4}$. 

We define a {\sl reduction} of a cutting sequence to be a replacement of the 
sequence by a shorter one, according to the one of the following rules
(where $\updownarrow$ denotes $\up$ {\it or} $\down$, and 
$i\in \{0,\ldots n\}$).

\begin{itemize}
\item[$\bullet$] $\underline{i}\updownarrow i \ \to \ \underline{i}$, \ 
$\underline{i+1}\updownarrow i \ \to \ \underline{i+1}$, \
$i\updownarrow \underline{i} \ \to \ \underline{i}$, \
$i\updownarrow \underline{i+1} \ \to \ \underline{i+1}$, 
\item[$\bullet$] $\down i \down \ \to \ \down$, \ $\up i \up \ \to \ \up$,
\item[$\bullet$] $i \updownarrow i \ \to \ i$, 
\item[$\bullet$] $\underline{i} \updownarrow \underline{i+1} \ \to \ 
\underline{i}\ \underline{i+1}$, \ $\underline{i+1} \updownarrow \underline{i} 
\ \to \ \underline{i+1}\ \underline{i}$ 
\end{itemize}

A cutting sequence is called {\sl reduced} if it allows no reduction.
\begin{prop}  Every braid on $n$ strings has a unique reduced cutting 
sequence. \end{prop}
\begin{proof} Let $\chi$ be a cutting sequence of a curve diagram $\Gamma$ of 
the braid.
We observe that a reduced version $\chi'$ of $\chi$ is the same as the cutting
sequence of a curve diagram $\Gamma'$, where $\Gamma'$ is obtained 
by reducing $\Gamma$ \wrt the trivial curve diagram $E$. From proposition
\ref{welldef} we deduce that any two reduced cutting sequences $\chi'$ and 
$\chi''$
must come from curve diagrams which are equivalent \wrt $E$. Therefore
$\chi'$ and $\chi''$ must agree. \end{proof}

The reduced curve diagram can be reconstructed from
the reduced cutting sequence.  Thus the cutting sequence classifies
the curve diagram, and hence the braid. This is most easily seen by using
pen and paper. One reads the cutting sequence, and for every
number symbol one encounters, draws one arc in the diagram.
If the cutting sequence is reduced, then this involves no choices. 
Below we shall give an algorithm to do this which is more suitable for 
computer implementation.  

Note that it is easy to construct reduced cutting 
sequences which do not come from curve diagrams.  The pen and paper 
method can also be used to decide whether a cutting sequence does 
correspond to a curve diagram.  Again we give a more formal algorithm 
below which will do this.

\sh{Reading the cutting sequence from the braid word}

We next show how to convert a braid defined in terms of the twist 
generators $\sigma_i^{\pm 1}$ into a reduced cutting sequence. We do 
this inductively by defining how $\sigma_i$ and $\sigma_i\inv$ act on 
reduced cutting sequences
and then let the whole word act on the trivial sequence 
$\underline{0}\ \underline{1}\ldots \underline{n}\ \underline{n+1}$.

\begin{algo}\label{braidtocut}
Suppose a braid $\beta$ has reduced cutting sequence $\chi$. Then a cutting
sequence of $\beta \sigma_i$ is obtained by simultaneously making the 
following replacements everywhere in the word $\chi$. These rules are to
be interpreted as simultaneous, not sequential, replacements.
\end{algo}

\begin{itemize}
\item[(i)] $\underline{i}\ \to \ \underline{i+1}$, \ \ $\underline{i+1}
\ \to \ \underline{i}$,
\item[(ii)] $\down (\underline{i})\ \to \ \down i-1 \up (\underline{i+1})$, 
\ \ $(\underline{i}) \down \ \to \ (\underline{i+1}) \up i-1 \down$,
\item[(iii)] $\underline{i-1} \ (\underline{i}) \ \to \ 
\underline{i-1} \up (\underline{i+1})$, \ \ $(\underline{i}) \ 
\underline{i-1} \ \to \ (\underline{i+1}) \up \underline{i-1}$,
\item[(iv)] $\up (\underline{i}) \ \to \ \up (\underline{i+1})$, \ \
$(\underline{i}) \up \ \to \  (\underline{i+1}) \up$,
\item[(v)] $\down (\underline{i+1})\ \to \ \down (\underline{i})$, \ \
$(\underline{i+1)} \down \ \to \ (\underline{i}) \down$,
\item[(vi)] $\underline{i+2} \ (\underline{i+1}) \ \to \ \underline{i+2} 
\down (\underline{i})$, \ \ $(\underline{i+1}) \ \underline{i+2} \ 
\to \ (\underline{i}) \down \underline{i+2}$,
\item[(vii)] $\up (\underline{i+1}) \ \to \ \up i+1 \down (\underline{i})$, 
\ \ $(\underline{i+1}) \up  \ \to \ (\underline{i}) \down i+1 \up$.
\item[(viii)] $\down i \up \ \to \ \down i-1 \up i \down i+1 \up$, \ \
$\up i \down \ \to \ \up i+1 \down i \up i-1 \down$,
\end{itemize}

{\sl Note: in rules {\rm (ii) - (vii)}, rule {\rm (i)} is being applied,
and its application is indicated by brackets. Replacements of symbols
other than $\underline{i}, \underline{i+1}$ depend on context, eg rule
{\rm (ii)} says that if $\down$ is followed by $\underline{i}$, 
then it is to be replaced by $\down i-1 \up$, and the $\underline{i}$ is 
replaced by $\underline{i+1}$, by {\rm (i)}. So $\down \underline{i}$ 
turns into $\down i-1 \up \underline{i+1}$.

The rules for the action of $\sigma_i\inv$ are obtained by interchanging
the symbols $\up$ and $\down$ everywhere in this list (ie replacing
up- by down-, and down- by up-arrows).
The resulting cutting sequence can then be reduced, to obtain the reduced 
cutting sequence of the braid $\beta \sigma_i$ or $\beta \sigma_i\inv$.} 
\ppar

We can now deduce an effective algorithm to decide whether a given braid
is positive, trivial, or negative:

\begin{algo}[To decide if a given braid is positive, trivial, or 
negative]\nl Use algorithm \ref{braidtocut} to
calculate the reduced cutting sequence of the braid. The braid is
positive if and only if the first arrow in this sequence is an
up-arrow $\up$.
\end{algo}

{\bf Recovering the curve diagram from the cutting sequence}
\ppar

We now show how to recover a reduced curve diagram from its
associated cutting sequence.  At the same time this will provide
an effective algorithm to decide if a given cutting sequence 
corresponds to a curve diagram.  

To make precise the problem here, we define the {\sl real cutting 
sequence}
of a curve diagram to be the cutting sequence, with the non-underlined
integers replaced by real numbers specifying the precise intersection point
of the curve diagram with the real line, up to order preserving bijections.
(Taking the integer part of all numbers in the real cutting sequence we
retrieve the cutting sequence.) Given the real cutting sequence, we can 
immediately construct the curve diagram.  Moreover it is trivial to
check if a real cutting sequence corresponds to an (embedded) curve
diagram:  one just checks that 

(1)\stdspace if $\underline{i}\ \underline{i+1}$ or 
$\underline{i+1}\ \underline{i}$ occurs in the sequence then
no real number in $(i,i+1)$ occurs,

(2)\stdspace the numbers on each side of two arrows
of the same type correspond to nested intervals (so that the corresponding
curves do not intersect).

So we need an algorithm to 
reconstruct the real cutting sequence from the cutting sequence
or equivalently to decide for each $i$ the order in which the
corresponding points actually occur in $\R$.  

\begin{algo} \label{reconstruct}
Suppose the letter $i$ $(i\in \{0,\ldots, n\})$ appears in two different 
places, say in the $r$th and $s$th position, in the cutting sequence.
To decide which one represents the smaller 
number in the interval $(i,i+1)$ in the real cutting sequence
proceed as follows. 
\end{algo}
%\vspace{-4mm}

Since the cutting sequence is reduced,
there are two arrows in opposite direction adjacent to each of the letters
$i$. Starting at the $r$th letter we read the sequence either
forwards or backwards. We define the {\sl up-string at the $r$th place} to 
be the word obtained from the cutting sequence by reading forwards or 
backwards, starting at the $r$th letter, up to the next underlined number,
with the reading direction specified by the requirement that the 
the first two letters read should be $i \up$. Similarly, we define
the {\sl down-string at the $r$th place} by reading in the opposite direction,
such that the resulting word starts with $i \down$, again up to the next
underlined number. We compare the up-string at the $r$th with that at
the $s$th place, and the down-string at the $r$th with that at the
$s$th place. They cannot both agree, for if they did, the curve diagram
would have two curves with the same endpoints.

We now manipulate the up- and down strings as follows: firstly, we increase
all non-underlined integers by $\half$. Then we remove the underline from
all underlined integers. We obtain sequences of the form
$x_0 \updownarrow x_1 \updownarrow \ldots \updownarrow x_{l-1} \updownarrow
x_l$, where $l\in \N$, $x_0=i+\half$, $x_1,\ldots, x_{l-1} \in 
\{\half, 1\half, \ldots, n+\half\}$, and $x_l \in \{0,\ldots n+1\}$.

 From this we can construct  a sequence of numbers in $\{1,1\half,\ldots,
n-\half, n\}$, called the {\sl cyclically associated sequence}, as follows.
For every string $x_j \up x_{j+1}$ we write down the unique representative
in $\{\half, 1,\ldots,n, n+\half\}$ of $x_{j+1} - x_j + (n+1)\Z
\in \R / (n+1)\Z$;
for every string $x_j \down x_{j+1}$ we write down the unique representative
in $\{\half,1 ,\ldots,n, n+\half\}$ of $x_j - x_{j+1} + (n+1)\Z \in 
\R / (n+1)\Z$. Altogether, this yields a sequence of length $l$.

We now define an up-string $u$ to be {\sl cyclically
lexicographically larger} than another up-string $u'$, if the cyclically 
associated sequence of $u$ is lexicographically larger than the one of 
$u'$.\fnote{Cyclic lexicographic order is used by Birman and Series 
\cite{BirSer}.}
The geometric interpretation is that the curve diagram has two line segments 
starting in the real interval $(i,i+1)$, going into the upper half plane.
The line segment representing the cyclically lexicographically larger
up-string is the one turning `more to the left'. Since the two line
segments must be disjoint (being part of the curve diagram), the starting
point of the curve segment yielding the cyclically lexicographically
larger up-string must represent a smaller real number in the real cutting
sequence. Similarly, we define a cyclic lexicographic ordering on the 
down-strings; this time, the starting point of a curve segment which 
gives rise to a cyclically lexicographically larger down-string than another 
curve segment must represent a {\it larger} real number in the real cutting
sequence.\stdspace {\sl End of algorithm \ref{reconstruct}}

To summarise, we have found an algorithm for reconstructing the real cutting
sequence from the cutting sequence: given any two places in the cutting
sequence where the letter $i$ occurs, we compare the up-strings at these
places. If they agree, we compare the down-strings instead. In either
case we can work out the cyclically associated sequences, and then
decide which of the two letters $i$ represents the smaller number in the
interval $(i,i+1)$ in the real cutting sequence. 

\sh{An algorithm to determine order from the cutting sequence}

Algorithm \ref{reconstruct} also allows us to decide which of two given
reduced cutting sequences represents the larger braid. If the two sequences 
agree on some initial segment, then we remove the underlines from all
underlined numbers (except the first letter $\underline{0}$) that lie in 
this segment. Then we reduce the resulting two sequences. We obtain two new 
sequences whose initial segments up to the first underlined numbers
do not agree. If they differ already on the second letter (after 
$\underline{0}$), then we know which one is larger. Otherwise, we work out 
which of them is cyclically lexicographically larger, using algorithm 
\ref{reconstruct}. 

\sh{The algorithm to determine left-consistent canonical form}

We are finally ready to describe our algorithm to calculate the 
left-consistent canonical form of a braid.  The input is a braid
$\beta$ represented as a word $w$ in the twist generators 
$\sigma_i^{\pm 1}$.  The output is the same braid in
left-consistent canonical form of $\beta$, again given as a word
in the $\sigma_i^{\pm 1}$.

The algorithm proceeds by repeating the main step (described below)
after each repetition we have a word $W$ and a cyclically reduced
cutting sequence $\chi$ which are both modified at the next
repetition.

\rk{Start}We start with $W$ the trivial word, and $\chi$ the
reduced cutting sequence of $\beta$ calculated using algorithm
\ref{braidtocut}.

\rk{Finish}If the reduced cutting sequence $\chi$ is 
$\underline{0}\ \underline{1}\ \ldots \underline{n}\ \underline{n+1}$, 
then the algorithm stops, and the inverse of the word $W$ is the desired 
canonical word. 

\rk{Main step}
If the reduced cutting sequence starts $\underline{0}\ \underline{1} \ldots
\underline{i} \up$, with $i<n+1$, then we hunt for subwords of the 
following forms

\begin{itemize}
\item[(i)] $i \up a_1 \down a_2 \up \ldots \updownarrow a_{l-1} \updownarrow
\underline{a_l}$ \ or
\item[(ii)] $\underline{i} \up a_1 \down a_2 \up \ldots \updownarrow a_{l-1} 
\updownarrow \underline{a_l}$ \ or
\item[(iii)] $\underline{a_l} \updownarrow a_{l-1} \updownarrow \ldots \up 
a_2 \down a_1 \up i$ \ or
\item[(iv)] $\underline{a_l} \updownarrow a_{l-1} \updownarrow \ldots \up 
a_2 \down a_1 \up \underline{i}$,
\end{itemize}

where the $a_1,\ldots, a_{l-1}$ are not equal to $i$ and not underlined,
and $a_l \ne i,i+1$. (If the reduced cutting sequence starts $\underline{0}
\ldots \underline{i} \down$, then we hunt for subwords like $i\down a_1
\up \ldots \updownarrow a_{l-1} \updownarrow a_l$ instead.)
We shall call these words {\sl useful subwords}, because they correspond to
useful arcs. 

We consider the set of all useful subwords, and we want to
identify the `leftmost one', ie the one whose letter $i$ represents the
leftmost point in the interval $(i,i+1)$. If one of them starts or ends
with a letter $\underline{i}$, ie if one of them is of type (ii) or (iv),
then this is it. If not, then we can use algorithm \ref{reconstruct} to determine the leftmost one.
When we have found the leftmost useful subword, we modify it as follows.
If it is of type (i) or (ii), then we write it backwards, so that it starts
with the letter $\underline{a_l}$. Irrespectively of the type of the
useful subword, we remove the underline from the letter $\underline{a_l}$.
Then we let $c:=a_l$, replace all letters $a_k$ ($k\in \{1,\ldots, l\}$) 
with $a_k \geqslant c$ by $a_k-1$ (eg $a_l$ turns into $a_l-1$), and reduce
the resulting sequence. By doing this, we obtain a modified sequence 
$a'_0 \updownarrow a'_1 \updownarrow \ldots \updownarrow a'_{l'}$, possibly 
with the letter $a'_{l'}=i$ underlined.

We now multiply $W$ on the right by a word $v_1 \ldots v_l$, where $v_k$
is determined by $a'_{k-1}$, $a'_k$, and the arrow in between $a'_{k-1}$
and $a'_k$ as follows:

\begin{itemize}
\item[(i)] If the modified leftmost useful subword contains the string
$a'_{k-1} \up a'_k$, and $a'_{k-1} < a'_k$, then 
$v_k=\sigma_{a'_{k-1}+1} \ldots \sigma_{a'_k}$;
\item[(ii)] If the modified leftmost useful subword contains the string
$a'_{k-1} \up a'_k$, and $a'_{k-1} > a'_k$, then 
$v_k=\sigma_{a'_{k-1}}\inv \ldots \sigma_{a'_k+1}\inv$;
\item[(iii)] If the modified leftmost useful subword contains the string
$a'_{k-1} \down a'_k$, and $a'_{k-1} < a'_k$, then 
$v_k=\sigma_{a'_{k-1}+1}\inv \ldots \sigma_{a'_k}\inv$;
\item[(iv)] If the modified leftmost useful subword contains the string
$a'_{k-1} \down a'_k$, and $a'_{k-1} > a'_k$, then 
$v_k=\sigma_{a'_{k-1}} \ldots \sigma_{a'_k+1}$
\end{itemize}

The word $v_1 \ldots v_l$ represents the slide of a hole back along the 
leftmost useful arc.

Finally, we calculate the new reduced cutting sequence after this slide.
This can be done by letting the word $v_1 \ldots v_l$ act on the reduced 
cutting sequence, as described above. (An alternative method would be to
remove the underline from the letter $\underline{a_l}$, underline the
unique letter $i$ which belongs to the leftmost useful subword instead,
carefully relabel the cutting sequence, using algorithm \ref{reconstruct}, 
and then reduce the resulting cutting sequence.)\stdspace
{\sl End of main step}

The proof of theorem \ref{lcform} implies that the algorithm stops
after a finite number of repetitions of the main step.


\begin{thebibliography}{999}


\bibitem{Bergman}{\bf G M Bergman} {\it Right orderable groups that are
not locally indicable}, Pacific J Math 174 (1991) 243--248

\bibitem{Birman}{\bf J  Birman}, {\it Braids, links, and mapping 
class groups}, Annals of Math. Studies, 82, Princeton University Press, 
Princeton (1975)

\bibitem{BirSer}{\bf J\,S Birman}, {\bf C Series}, {\it An algorithm
for simple curves on surfaces}, J. London Math. Soc (2) 29 (1984)
331--342

\bibitem{BuHa}{\bf R G Burns, V W D Hale}, {\it A note on group rings of
certain torsion free groups}, Canad Math Bull 15 (1972) 441--445

\bibitem{D.proofnoaxiom}{\bf P Dehornoy}, {\it Braid groups and left 
distributive operations}, Trans. AMS 345 (1994) 115--150

\bibitem{D.axiom}{\bf P Dehornoy}, {\it From large cardinals to braids 
via distributive algebra}, J. Knot Theory and its Ramifications 4(1995)
33--79

\bibitem{D.algorithm}{\bf P Dehornoy}, {\it A fast method of comparing 
braids}, Adv. in Math. 125 (1997) 200--235

\bibitem{elfir}{\bf E\,A Elrifai}, {\bf H\,R Morton}, {\it Algorithms 
for positive braids},
Quart. J. Math. Oxford 45  (1994) 479--497

\bibitem{Ep}{\bf D\,B\,A Epstein et al}, {\it Word processing in groups},
Jones \& Bartlett (1992)

\bibitem{GoLin}{\bf E A Gorin, V Ja Lin} {\it Algebraic equations with 
continuous coefficients, and certain questions of the algebraic theory of 
braids}, Math USSR-Sb 7 (1969) 569-596

\bibitem{laver}{\bf R Laver}, {\it Braid group actions on left-distibutive
structures and well-orderings in the braid group}, J. Pure Appl. Algebra
108 (1996) 81--98

\bibitem{Mosher1} {\bf L Mosher}, {\it Mapping class groups are 
automatic}, Math. Research Letters 1 (1994) 249--255

\bibitem{Mosher2} {\bf L Mosher}, {\it Mapping class groups are 
automatic}, Annals of Math. 142 (1995) 303--384

\bibitem{Mosher.practical} {\bf L Mosher}, {\it A user's guide to the
mapping class group: once punctured surfaces}, MSRI preprint


\bibitem{lrinvariant}{\bf D Rolfsen}, {\bf Jun Zhu}, {\it Braids, 
orderings and zero divisors},  submitted to J. Knot Theory and
its Ramifications

\bibitem{Rourke-Wiest}{\bf C Rourke}, {\bf B Wiest}, {\it Order
automatic mapping class groups}, (to appear), 
{\tt http://www.maths.warwick.ac.uk/\char'176cpr/ftp/ordaut.ps}

\bibitem{bert}{\bf B Wiest}, {\it Dehornoy's ordering of the braid 
groups extends the subword ordering}, Pacific J. Math. (to appear) 

\end{thebibliography}
\end{document}